# DISCUSSION OF "LEAST ANGLE REGRESSION" BY EFRON ET AL.

By Saharon Rosset and Ji Zhu

*IBM T. J. Watson Research Center and Stanford University*

**1. Introduction.** We congratulate the authors on their excellent work. The paper combines elegant theory and useful practical results in an intriguing manner. The LAR–Lasso–boosting relationship opens the door for new insights on existing methods' underlying statistical mechanisms and for the development of new and promising methodology. Two issues in particular have captured our attention, as their implications go beyond the squared error loss case presented in this paper, into wider statistical domains: robust fitting, classification, machine learning and more. We concentrate our discussion on these two results and their extensions.

**2. Piecewise linear regularized solution paths.** The first issue is the piecewise linear solution paths to regularized optimization problems. As the discussion paper shows, the path of optimal solutions to the "Lasso" regularized optimization problem

$$(2.1) \qquad \hat{\beta}(\lambda) = \arg\min_\beta \|y - X\beta\|_2^2 + \lambda\|\beta\|_1$$

is piecewise linear as a function of $\lambda$; that is, there exist $\infty > \lambda_0 > \lambda_1 > \cdots > \lambda_m = 0$ such that $\forall \lambda \geq 0$, with $\lambda_k \geq \lambda \geq \lambda_{k+1}$, we have

$$\hat{\beta}(\lambda) = \hat{\beta}(\lambda_k) - (\lambda - \lambda_k)\gamma_k.$$

In the discussion paper's terms, $\gamma_k$ is the "LAR" direction for the $k$th step of the LAR–Lasso algorithm.

This property allows the LAR–Lasso algorithm to generate the whole path of Lasso solutions, $\hat{\beta}(\lambda)$, for "practically" the cost of one least squares calculation on the data (this is exactly the case for LAR but not for LAR–Lasso, which may be significantly more computationally intensive on some data sets). The important practical consequence is that it is not necessary to select the regularization parameter $\lambda$ in advance, and it is now computationally feasible to optimize it based on cross-validation (or approximate $C_p$, as presented in the discussion paper).

The question we ask is: what makes the solution paths piecewise linear? Is it the use of squared error loss? Or the Lasso penalty? The answer is that both play an important role. However, the family of (loss, penalty) pairs which facilitates piecewise linear solution paths turns out to contain many other interesting and useful optimization problems.

We now briefly review our results, presented in detail in Rosset and Zhu (2004). Consider the general regularized optimization problem

$$(2.2) \qquad \hat{\beta}(\lambda) = \arg\min_{\beta} \sum_i L(y_i, \mathbf{x}_i^t \beta) + \lambda J(\beta),$$

where we only assume that the loss $L$ and the penalty $J$ are both convex functions of $\beta$ for any sample $\{\mathbf{x}_i^t, y_i\}_{i=1}^n$. For our discussion, the data sample is assumed fixed, and so we will use the notation $L(\beta)$, where the dependence on the data is implicit.

Notice that piecewise linearity is equivalent to requiring that

$$\frac{\partial \hat{\beta}(\lambda)}{\partial \lambda} \in \mathcal{R}^p$$

is piecewise constant as a function of $\lambda$. If $L$, $J$ are twice differentiable functions of $\beta$, then it is easy to derive that

$$(2.3) \qquad \frac{\partial \hat{\beta}(\lambda)}{\partial \lambda} = -(\nabla^2 L(\hat{\beta}(\lambda)) + \lambda \nabla^2 J(\hat{\beta}(\lambda)))^{-1} \nabla J(\hat{\beta}(\lambda)).$$

With a little more work we extend this result to "almost twice differentiable" loss and penalty functions (i.e., twice differentiable everywhere except at a finite number of points), which leads us to the following *sufficient conditions for piecewise linear* $\hat{\beta}(\lambda)$:

1. $\nabla^2 L(\hat{\beta}(\lambda))$ is piecewise constant as a function of $\lambda$. This condition is met if $L$ is a piecewise-quadratic function of $\beta$. This class includes the squared error loss of the Lasso, but also absolute loss and combinations of the two, such as Huber's loss.
2. $\nabla J(\hat{\beta}(\lambda))$ is piecewise constant as a function of $\lambda$. This condition is met if $J$ is a piecewise-linear function of $\beta$. This class includes the $l_1$ penalty of the Lasso, but also the $l_\infty$ norm penalty.

2.1. *Examples.* Our first example is the "Huberized" Lasso; that is, we use the loss

$$(2.4) \qquad L(y, \mathbf{x}\beta) = \begin{cases} (y - \mathbf{x}^t \beta)^2, & \text{if } |y - \mathbf{x}^t \beta| \leq \delta, \\ \delta^2 + 2\delta(|y - \mathbf{x}\beta| - \delta), & \text{otherwise,} \end{cases}$$

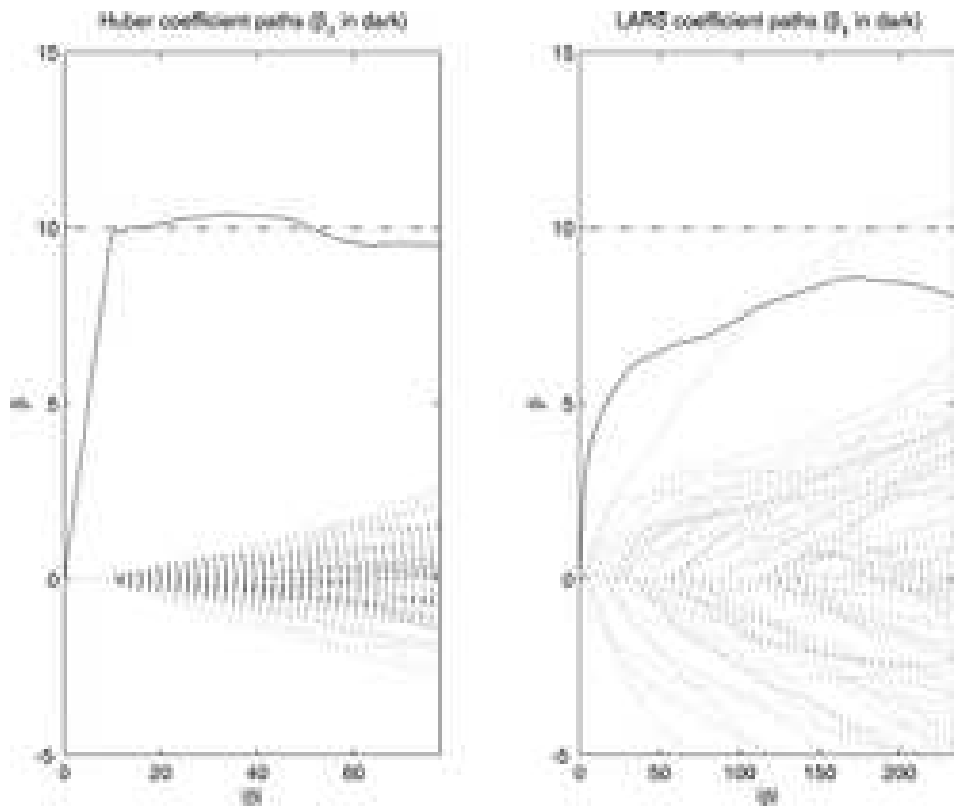

FIG. 1. *Coefficient paths for Huberized Lasso* (left) *and Lasso* (right) *for data example:* $\hat{\beta}_1(\lambda)$ *is the full line; the true model is* $E(Y|x) = 10x_1$.

with the Lasso penalty. This loss is more robust than squared error loss against outliers and long-tailed residual distributions, while still allowing us to calculate the whole path of regularized solutions efficiently.

To illustrate the importance of robustness in addition to regularization, consider the following simple simulated example: take $n = 100$ observations and $p = 80$ predictors, where all $x_{ij}$ are i.i.d. $N(0,1)$ and the true model is

$$(2.5) \qquad y_i = 10 \cdot x_{i1} + \varepsilon_i,$$

$$(2.6) \qquad \varepsilon_i \stackrel{\text{i.i.d.}}{\sim} 0.9 \cdot N(0,1) + 0.1 \cdot N(0,100).$$

So the normality of residuals, implicitly assumed by using squared error loss, is violated.

Figure 1 shows the optimal coefficient paths $\hat{\beta}(\lambda)$ for the Lasso (right) and "Huberized" Lasso, using $\delta = 1$ (left). We observe that the Lasso fails in identifying the correct model $E(Y|x) = 10x_1$ while the robust loss identifies it almost exactly, *if we choose the appropriate regularization parameter.*

As a second example, consider a classification scenario where the loss we use depends on the margin $y\mathbf{x}^t\beta$ rather than on the residual. In particular, consider the 1-*norm* support vector machine regularized optimization problem, popular in the machine learning community. It consists of minimizing the "hinge loss" with a Lasso penalty:

$$(2.7) \qquad L(y, \mathbf{x}^t\beta) = \begin{cases} (1 - y\mathbf{x}^t\beta), & \text{if } y\mathbf{x}^t\beta \leq 1, \\ 0, & \text{otherwise.} \end{cases}$$

This problem obeys our conditions for piecewise linearity, and so we can generate all regularized solutions for this fitting problem efficiently. This is particularly advantageous in high-dimensional machine learning problems, where regularization is critical, and it is usually not clear in advance what a good regularization parameter would be. A detailed discussion of the computational and methodological aspects of this example appears in Zhu, Rosset, Hastie, and Tibshirani (2004).

**3. Relationship between "boosting" algorithms and $l_1$-regularized fitting.** The discussion paper establishes the close relationship between $\varepsilon$-stagewise linear regression and the Lasso. Figure 1 in that paper illustrates the near-equivalence in the solution paths generated by the two methods, and Theorem 2 formally states a related result. It should be noted, however, that their theorem falls short of proving the global relation between the methods, which the examples suggest.

In Rosset, Zhu and Hastie (2003) we demonstrate that this relation between the path of $l_1$-regularized optimal solutions [which we have denoted above by $\hat{\beta}(\lambda)$] and the path of "generalized" $\varepsilon$-stagewise (AKA boosting) solutions extends beyond squared error loss and in fact applies to any convex differentiable loss.

More concretely, consider the following generic gradient-based "$\varepsilon$-boosting" algorithm [we follow Friedman (2001) and Mason, Baxter, Bartlett and Frean (2000) in this view of boosting], which iteratively builds the solution path $\beta^{(t)}$:

ALGORITHM 1 (Generic gradient-based boosting algorithm).

1. Set $\beta^{(0)} = 0$.
2. For $t = 1 : T$,

    (a) Let $j_t = \arg\max_j |\frac{\partial \sum_i L(y_i, \mathbf{x}_i^t \beta^{(t-1)})}{\partial \beta_j^{(t-1)}}|$.

    (b) Set $\beta_{j_t}^{(t)} = \beta_{j_t}^{(t-1)} - \varepsilon \operatorname{sign}(\frac{\partial \sum_i L(y_i, \mathbf{x}_i^t \beta^{(t-1)})}{\partial \beta_{j_t}^{(t-1)}})$ and $\beta_k^{(t)} = \beta_k^{(t-1)}$, $k \neq j_t$.

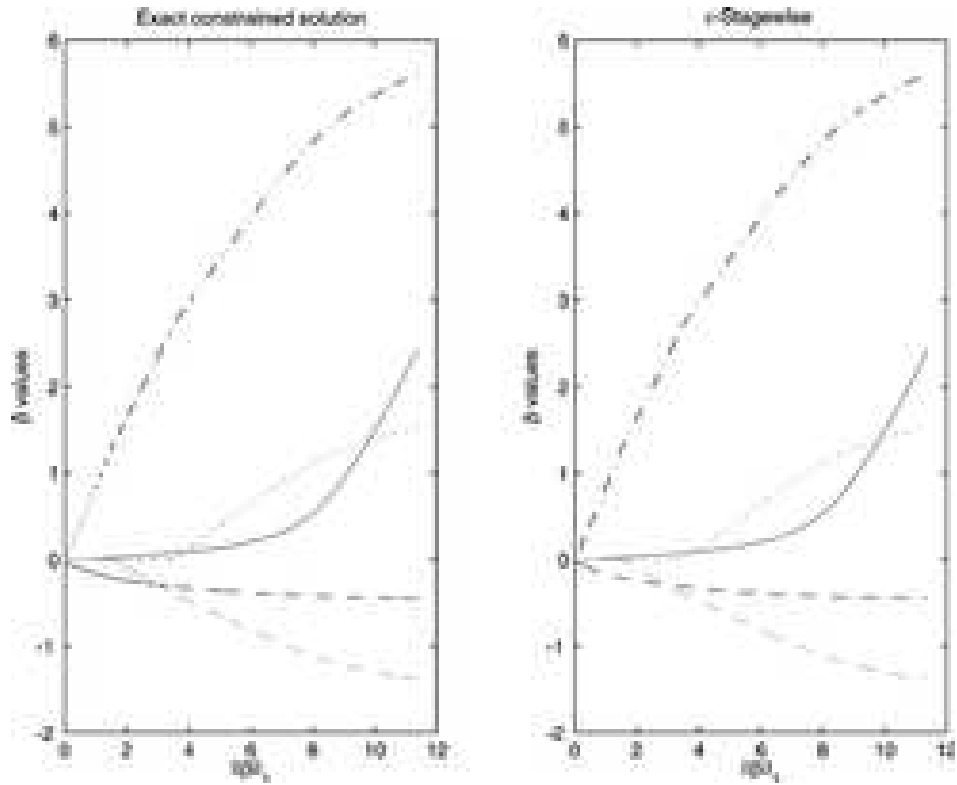

Fig. 2. *Exact coefficient paths* (left) *for $l_1$-constrained logistic regression and boosting coefficient paths* (right) *with binomial log-likelihood loss on five variables from the "spam" dataset. The boosting path was generated using $\varepsilon = 0.003$ and 7000 iterations.*

This is a coordinate descent algorithm, which reduces to forward stagewise, as defined in the discussion paper, if we take the loss to be squared error loss: $L(y_i, \mathbf{x}_i^t \beta^{(t-1)}) = (y_i - \mathbf{x}_i^t \beta^{(t-1)})^2$. If we take the loss to be the exponential loss,

$$L(y_i, \mathbf{x}_i^t \beta^{(t-1)}) = \exp(-y_i \mathbf{x}_i^t \beta^{(t-1)}),$$

we get a variant of AdaBoost [Freund and Schapire (1997)]—the original and most famous boosting algorithm.

Figure 2 illustrates the equivalence between Algorithm 1 and the optimal solution path for a simple logistic regression example, using five variables from the "spam" dataset. We can see that there is a perfect equivalence between the regularized solution path (left) and the "boosting" solution path (right).

In Rosset, Zhu and Hastie (2003) we formally state this equivalence, with the required conditions, as a conjecture. We also generalize the weaker result, proven by the discussion paper for the case of squared error loss, to any convex differentiable loss.

This result is interesting in the boosting context because it facilitates a view of boosting as approximate and implicit regularized optimization. The situations in which boosting is employed in practice are ones where explicitly solving regularized optimization problems is not practical (usually very high-dimensional predictor spaces). The approximate regularized optimization view which emerges from our results allows us to better understand boosting and its great empirical success [Breiman (1999)]. It also allows us to derive approximate convergence results for boosting.

**4. Conclusion.** The computational and theoretical results of the discussion paper shed new light on variable selection and regularization methods for linear regression. However, we believe that variants of these results are useful and applicable beyond that domain. We hope that the two extensions that we have presented convey this message successfully.

**Acknowledgment.** We thank Giles Hooker for useful comments.

IBM T. J. WATSON RESEARCH CENTER
P.O. BOX 218
YORKTOWN HEIGHTS, NEW YORK 10598
USA
E-MAIL: srosset@us.ibm.com

DEPARTMENT OF STATISTICS
UNIVERSITY OF MICHIGAN
550 EAST UNIVERSITY
ANN ARBOR, MICHIGAN 48109-1092
USA
E-MAIL: jizhu@umich.edu